\documentclass[11pt,reqno]{amsart}
\usepackage{amsmath,amsthm,amssymb,verbatim,enumerate,ifthen}
\usepackage[mathscr]{eucal}
\usepackage[utf8]{inputenc}
\def\N{\mathbb{N}}
\def\R{\mathbb{R}}

\def\I{\mathbb{I}}

\newtheorem{theorem}{Theorem}
\newtheorem*{theorem*}{Theorem}
\def\Thm#1#2{\ifthenelse{\equal{#1}{*}}{\begin{theorem*}#2\end{theorem*}}
  {\begin{theorem}\label{T#1}#2\end{theorem}}}
\newtheorem{Atheorem}{Theorem}

\def\thm#1{Theorem~\ref{T#1}}

\newtheorem{proposition}[theorem]{Proposition}
\newtheorem*{proposition*}{Proposition}
\def\Prp#1#2{\ifthenelse{\equal{#1}{*}}{\begin{proposition*}#2\end{proposition*}}
             {\begin{proposition}\label{P#1}#2\end{proposition}}}

\newtheorem{corollary}[theorem]{Corollary}
\newtheorem*{corollary*}{Corollary}
\def\Cor#1#2{\ifthenelse{\equal{#1}{*}}{\begin{corollary*}#2\end{corollary*}}
             {\begin{corollary}\label{C#1}#2\end{corollary}}}
\def\cor#1{Corollary~\ref{C#1}}

\newtheorem{lemma}[theorem]{Lemma}
\newtheorem*{lemma*}{Lemma}
\def\Lem#1#2{\ifthenelse{\equal{#1}{*}}{\begin{lemma*}#2\end{lemma*}}
             {\begin{lemma}\label{L#1}#2\end{lemma}}}
\def\lem#1{Lemma~\ref{L#1}}

\newtheorem{example}[theorem]{Example}
\newtheorem*{example*}{Example}
\def\Exa#1#2{\ifthenelse{\equal{#1}{*}}{\begin{example*}\rm #2\end{example*}}
             {\begin{example}\label{Ex#1}\rm #2\end{example}}}

\newtheorem{problem}[theorem]{Problem}

\theoremstyle{definition}
\newtheorem{definition}[theorem]{Definition}

\newtheorem{remark}[theorem]{Remark}
\newtheorem*{remark*}{Remark}
\def\Rem#1#2{\ifthenelse{\equal{#1}{*}}{\begin{remark*}\rm #2\end{remark*}}
             {\begin{remark}\label{R#1}\rm #2\end{remark}}}

\def\eq#1{{\rm(\ref{E#1})}}
\def\Eq#1#2{\ifthenelse{\equal{#1}{*}}
  {\begin{equation*}\begin{aligned}#2\end{aligned}\end{equation*}}
  {\begin{equation}\begin{aligned}\label{E#1}#2\end{aligned}\end{equation}}}

\long\def\comment#1{}

\begin{document}
\begin{flushright}
\end{flushright}
\vspace{5mm}

\date{\today}
\title{An extension of the Abel--Liouville identity}

\author[Zs. P\'ales]{Zsolt P\'ales}
\address[Zs. P\'ales]{Institute of Mathematics, University of Debrecen, H-4002 Debrecen, Pf.\ 400, Hungary}
\email{pales@science.unideb.hu}

\author[A. Zakaria]{Amr Zakaria}
\address[A. Zakaria]{Doctoral School of Mathematical and Computational Sciences, University of Debrecen, H-4002 Debrecen, Pf.\ 400, Hungary; Department of Mathematics, Faculty of Education, Ain Shams University, Cairo 11341, Egypt}
\email{amr.zakaria@science.unideb.hu}



\thanks{The research of the first author was supported by the K-134191 NKFIH Grant and the 2019-2.1.11-TÉT-2019-00049, the EFOP-3.6.1-16-2016-00022 and the EFOP-3.6.2-16-2017-00015 projects. The last two projects are co-financed by the European Union and the European Social Fund.}

\begin{abstract}
In this note, we present an extension of the celebrated Abel--Liouville identity in terms of noncommutative complete Bell polynomials for generalized Wronskians. We also characterize the range equivalence of $n$-dimensional vector-valued functions in the subclass of $n$-times differentiable functions with a nonvanishing Wronskian.
\end{abstract}

\maketitle

\section{\bf Introduction}

Throughout this paper let $\N$ and $\N_0$ denote the set of positive and nonnegative integers, respectively, and let $I$ stand for a nonempty open real interval.

For an $n$-dimensional vector valued $(n-1)$-times continuously differentiable function $f:I\to\R^n$, its \emph{Wronskian} $W_f:I\to\R$ is defined by 
\Eq{*}{
  W_f:=\begin{vmatrix}
  f^{(n-1)} & \dots & f' & f
  \end{vmatrix}.
}
Consider now the $n$th-order homogeneous linear differential equation 
\Eq{DE}{
  y^{(n)}
  =a_{1}y^{(n-1)}+\dots+a_{n}y,
}
where $a_1,\dots,a_n:I\to\R$ are continuous functions. By the classical \emph{Abel--Liouville identity} (cf.\ \cite{Tes12}), if $f:I\to\R^n$ is a fundamental system of solutions of \eq{DE}, then $W_f$ does not vanish on $I$ and 
\Eq{*}{
  W_f'=a_1W_f.
}

For a sufficiently smooth function $f:I\to\R^n$ and $k=(k_1,\dots,k_n)\in\N_0^n$, we introduce now the \emph{generalized Wronskian} $W_f^k:I\to\R$ by
\Eq{*}{
  W_f^k:=\begin{vmatrix}
  f^{(k_1)} & \dots & f^{(k_n)} 
  \end{vmatrix}.
}
One can easily see that, with this notation, we have
\Eq{*}{
  W_f=W_f^{(n-1,n-2,\dots,0)}\qquad \mbox{and}\qquad
  W_f'=W_f^{(n,n-2,\dots,0)}.
}
Therefore, the Abel--Liouville identity can be rewritten as
\Eq{L}{
  W_f^{(n,n-2,\dots,0)}=a_1W_f^{(n-1,n-2,\dots,0)}.
}
One of the main goal of this short paper is to establish a formula for $W_f^k$ in terms of the coefficients of differential equation \eq{DE}.
Another goal is to introduce the range equivalence for $n$-dimensional vector-valued functions and to characterize this equivalence relation in the subclass of $n$-times differentiable functions with a nonvanishing Wronskian.

\section{Main results}

For the description of our main result, we recall the notion of \emph{noncommutative complete Bell polynomials}, which was introduced by Schimming and Rida \cite{SchRid96}. Let $\R^{n\times n}$ denote the ring of $n\times n$ matrices with real entries and let $\I_n$ denote the $n\times n$ unit matrix. Now define $B_m:(\R^{n\times n})^m\to \R^{n\times n}$ by the following recursive formula
\Eq{Bell}{
  B_0:=\I_n,\qquad B_{m+1}(X_1,\dots,X_{m+1}):=\sum_{j=0}^m \binom{m}{j} B_j(X_1,\dots,X_j)X_{m+1-j}.
}
The notion of \emph{complete Bell polynomials} in the commutative setting (i.e., when $n=1$) was introduced by Bell \cite{Bel27}, \cite{Bel34}. One can easily compute the first few Bell polynomials as follows:
\Eq{*}{
 B_1(X_1)&=X_1,\\
 B_2(X_1,X_2)&=X_1^2+X_2,\\
 B_3(X_1,X_2,X_3)&=X_1^3+2X_1X_2+X_2X_1+X_3,\\
 B_4(X_1,X_2,X_3,X_4)&=X_1^4+3X_1^2X_2+2X_1X_2X_1+3X_1X_3+3X_2^2+X_2X_1^2+X_3X_1
 +X_4,\\
 B_5(X_1,X_2,X_3,X_4,X_5)
 &=X_1^5+4X_1^3X_2+3X_1^2X_2X_1+6X_1^2X_3+8X_1X_2^2+2X_1X_2X_1^2\\
 &\qquad+3X_1X_3X_1+4X_1X_4+3X_2^2X_1+X_2X_1^3+4X_2X_1X_2+6X_2X_3\\
 &\qquad+X_3X_1^2+4X_3X_2+X_4X_1+X_5.}

The statement of the next basic lemma was proved in the paper \cite{SchRid96}.

\Lem{0}{For every $j\in\N_0$, and $j$-times differentiable matrix-valued function $X:I\to\R^{n\times n}$, 
\Eq{00}{
  B_{j+1}\big(X,\dots,X^{(j)}\big)=
   XB_{j}\big(X,\dots,X^{(j-1)}\big)
   +\Big(B_{j}\big(X,\dots,X^{(j-1)}\big)\Big)'.
}}

\Lem{1}{Let $n,m\in\N$, let $X:I\to\R^{n\times n}$ be an $(m-1)$-times continuously differentiable function 
and $Y:I\to\R^{n\times n}$ be a  differentiable function such that 
\Eq{Y1}{
  Y'=YX
}
holds on $I$. Then $Y$ is $m$-times continuously differentiable and
\Eq{Yk}{
Y^{(j)}=Y B_j\big(X,\dots,X^{(j-1)}\big)
\qquad (j\in\{0,\dots,m\}).
}}

\begin{proof} If $m=1$, then $X$ is continuous, hence the continuity of $Y$ and \eq{Y1} imply that $Y$ is continuously differentiable. If $m>1$, then using \eq{Y1}, a simple inductive argument shows that $Y$ is $m$-times 
continuously differentiable.

The equality \eq{Yk} is trivial if $j=0$, because $B_0=\I_n$. For $j=1$, the equality \eq{Yk} is equivalent to \eq{Y1}. 
Now assume that \eq{Yk} has been verified for some $1\leq j<m$. Then, using \eq{Y1} and \lem{0}, we get
\Eq{*}{
  Y^{(j+1)}&=\big(Y^{(j)}\big)'
  =\Big(YB_{j}\big(X,\dots,X^{(j-1)}\big)\Big)'
  =Y'B_{j}\big(X,\dots,X^{(j-1)}\big)
     +Y\Big(B_{j}\big(X,\dots,X^{(j-1)}\big)\Big)'\\
  &=Y\Big[XB_{j}\big(X,\dots,X^{(j-1)}\big)
   +\Big(B_{j}\big(X,\dots,X^{(j-1)}\big)\Big)'\Big]
  =  Y B_{j+1}\big(X,\dots,X^{(j)}\big).
}
This proves the assertion for $j+1$.
\end{proof}

In what follows, let $e_1,\dots,e_n$ denote the elements of the standard basis in $\R^n$.

\Cor{0}{Let $n,m\in\N$, let $a:I\to\R^n$ be an $(m-1)$-times continuously differentiable function and let $f:I\to\R^n$ be 
a fundamental system of solutions of the differential equation \eq{DE}. Let the matrix-valued functions $X_a:I\to\R^{n\times n}$
and $Y_f:I\to\R^{n\times n}$ be defined by
\Eq{XY}{
X_a:= \begin{pmatrix}a&e_1&\dots&e_{n-1}\end{pmatrix}
\qquad\mbox{and}\qquad
Y_f:=
\begin{pmatrix} f^{(n-1)} & \dots & f' & f  \end{pmatrix}.
}
Then $Y_f$ is $m$-times continuously differentiable and
\Eq{Bk}{
Y_f^{(j)}=Y_f B_j\big(X_a,\dots,X_a^{(j-1)}\big)
\qquad (j\in\{0,\dots,m\}).
}}

\begin{proof} The function $f$ satisfies the differential equation \eq{DE}, therefore $f^{(n)}=Y_f\cdot a$. 
On the other hand, $f^{(n-i)}=Y_f\cdot e_i$ holds for $i\in\{1,\dots,n-1\}$. These equalities imply that 
\Eq{B1}{
Y_f'=\begin{pmatrix}f^{(n)} & f^{(n-1)} & \dots & f'\end{pmatrix}
=\begin{pmatrix}Y_f\cdot a & Y_f\cdot e_{1} &\dots & Y_f\cdot e_{n-1}\end{pmatrix}
=Y_f X_a.
}
Therefore, equation \eq{Y1} holds with $Y:=Y_f$ and $X:=X_a$, consequently, the statement is a consequence of \lem{1}.
\end{proof}

Using the above corollary, we can easily establish a formula for the computation of the generalized Wronskian $W_f^k$.

\Thm{1}{Let $n,m\in\N$, let $a:I\to\R^n$ be an $(m-1)$-times continuously differentiable function and let $f:I\to\R^n$ be 
a fundamental system of solutions of the differential equation \eq{DE}. Let the matrix-valued functions $X_a:I\to\R^{n\times n}$ 
be defined by \eq{XY}. Then, for $k=(k_1,\dots,k_n)\in\N_0^n$ with $\max(k_1,\dots,k_n)\leq m+n-1$, 
\Eq{Wk}{
  W_f^k
  =W_f\begin{vmatrix}
   B_{\ell_1}\big(X_a,\dots,X_a^{(\ell_1-1)}\big)e_{n+\ell_1-k_1} &\dots& 
   B_{\ell_n}\big(X_a,\dots,X_a^{(\ell_n-1)}\big)e_{n+\ell_n-k_n}
   \end{vmatrix},
}
where, for $i\in\{1,\dots,n\}$, $\ell_i:=(k_i-n+1)^+$.}

\begin{proof} Define the matrix valued function $Y_f:I\to\R^{n\times n}$ by \eq{XY} and observe that, by \cor{0}, 
for all $\ell\in\{0,\dots,m+n-1\}$, we have that
\Eq{*}{
   f^{(\ell)}
   =Y_f^{(i)}e_{n+i-\ell}
   =Y_f B_{i}\big(X_a,\dots,X_a^{(i-1)}\big)e_{n+i-\ell} 
   \qquad(i\in\{(\ell-n+1)^+,\dots,\min(\ell,m)\}).
}
By taking the smallest possible value for $i$ in the above formula, we get
\Eq{*}{
   f^{(\ell)}
   =Y_f B_{(\ell-n+1)^+}\big(X_a,\dots,X_a^{((\ell-n+1)^+-1)}\big)
     e_{n+(\ell-n+1)^+-\ell}.
}
Applying this equality for $\ell\in\{k_1,\dots,k_n\}$, we obtain
\Eq{*}{
  \begin{pmatrix} f^{(k_1)}&\dots&f^{(k_1)} \end{pmatrix}
  =Y_f\begin{pmatrix}
   B_{\ell_1}\big(X_a,\dots,X_a^{(\ell_1-1)}\big)e_{n+\ell_1-k_1} &\dots& 
   B_{\ell_n}\big(X_a,\dots,X_a^{(\ell_n-1)}\big)e_{n+\ell_n-k_n}
   \end{pmatrix}.
}
Now taking the determinant side by side and using the product rule for determinants, the equality \eq{Wk} follows. 
\end{proof}

In the subsequent corollary, we consider the case when $\ell_i=0$ for $i\in\{2,\dots,n\}$. In this particular setting, 
the determinant on the left hand side of \eq{Wk} can easily be computed.

\Cor{1}{Let $n,m\in\N$, let $a:I\to\R^n$ be an $(m-1)$-times continuously differentiable function and let $f:I\to\R^n$ be 
a fundamental system of solutions of the differential equation \eq{DE}. Let the matrix-valued functions $X_a:I\to\R^{n\times n}$
be defined by \eq{XY} and let $d\in\{0,...,m-1\}$ and $j\in\{0,\dots,n-1\}$. Then 
\Eq{dj}{
  W_f^{(n+d,n-1,\dots,j+1,j-1,\dots,0)}
  =(-1)^{n-j-1}W_f \langle B_{d+1}\big(X_a,\dots,X_a^{(d)}\big) e_1,e_{n-j}\rangle.
}
If $d=0$ and $j=n-1$, then this equality reduces to the Abel--Liouville identity \eq{L}. More generally, for $d=0,1,2$, we get the 
following formulas:
\Eq{12}{
  W_f^{(n,n-1,\dots,j+1,j-1,\dots,0)}
  &=(-1)^{n-j-1}W_f a_{n-j}, &\\
  W_f^{(n+1,n-1,\dots,j+1,j-1,\dots,0)}
  &=(-1)^{n-j-1}W_f (a_1a_{n-j}+a_{n-j+1}+a_{n-j}'), &\\
  W_f^{(n+2,n-1,\dots,j+1,j-1,\dots,0)}
  &=(-1)^{n-j-1}W_f (a_1^2a_{n-j}+a_1a_{n-j+1}+a_2a_{n-j}+a_{n-j+2}&\\ &\hspace{3cm} +a_1a_{n-j}'+2a_1'a_{n-j}+2a_{n-j+1}'+a_{n-j}'').
}
(Here we define $a_{n+1}:=a_{n+2}:=0$.)
}

\begin{proof}
We apply the previous theorem for $k:=(n+d,n-1,\dots,j+1,j-1,\dots,0)$, where $d\in\{0,...,m-1\}$ and $j\in\{0,\dots,n-1\}$. 
Then we get that $\ell_1=d+1$, and $\ell_i=0$ for $i\in\{2,\dots,n\}$. Therefore, 
\Eq{*}{
W_f^{(n+d,n-1,\dots,j+1,j-1,\dots,0)}
   &=W_f\begin{vmatrix}
   B_{d+1}\big(X_a,\dots,X_a^{(d)})e_1 &\I_ne_1&\dots&\I_ne_{n-j-1}&\I_ne_{n-j+1}&\dots& \I_ne_n
   \end{vmatrix}\\
   &=(-1)^{n-j-1}W_f \langle B_{d+1}\big(X_a,\dots,X_a^{(d)}\big) e_1,e_{n-j}\rangle.
}
Thus, equality \eq{dj} has been shown. In the case $d=0$, we have that 
\Eq{*}{
\langle B_1(X_a)e_1,e_{n-j}\rangle
=\langle X_ae_1,e_{n-j}\rangle=a_{n-j}
} 
because the $(n-j)$th entry of $X_a$ equals $a_{n-j}$. This implies the first equality in \eq{12} for $j\in\{0,\dots,n-1\}$. In particular, for $j=n-1$, this equality is equivalent to the Abel--Liouville identity \eq{L}.

In the case $d=1$, a simple computation gives that 
\Eq{*}{
\langle B_2(X_a,X_a')e_1,e_{n-j}\rangle
=\langle (X_a^2+X_a')e_1,e_{n-j}\rangle
=a_1a_{n-j}+a_{n-j+1}+a_{n-j}',
}
which yields the second equality in \eq{12} for $j\in\{0,\dots,n-1\}$.

In the case $d=2$, a somewhat more difficult computation gives that
\Eq{*}{
\langle B_3(X_a,X_a',X_a'')e_1,e_{n-j}\rangle
&=\langle (X_a^3+2X_aX_a'+X_a'X_a+X_a'')e_1,e_{n-j}\rangle \\
&=a_1^2a_{n-j}+a_1a_{n-j+1}+a_2a_{n-j}+a_{n-j+2}\\
&\hspace{3cm}+a_1a_{n-j}'+2a_1'a_{n-j}+2a_{n-j+1}'+a_{n-j}'',
}
which then yields the third equality in \eq{12}.
\end{proof}

For the sake of convenience and brevity, we introduce the following notation: 
for an $n$-times continuously differentiable function $f:I\to\R^n$ such that $W_f$ is nonvanishing and $j\in\{0,\dots,n-1\}$, the function $\Phi_f^{[j]}:I\to\R$ is defined by
\Eq{*}{
\Phi_f^{[j]}:=(-1)^{n-j-1}\frac{W_f^{(n,\dots,j+1,j-1,\dots,0)}}{W_f}.
}
For instance, if $f$ is $n$-times continuously differentiable function whose components form a fundamental system of solutions for \eq{DE}, then the Abel--Liouville identity \eq{L} can be rewritten as
\Eq{*}{
\Phi_f^{[n-1]}=a_1.
}
More generally, the first equality in \eq{12} gives that
\Eq{*}{
 \Phi_f^{[j]}=a_{n-j} 
 \qquad (j\in\{0,\dots,n-1\})
}
or, equivalently, 
\Eq{aj}{
 a_j=\Phi_f^{[n-j]} 
 \qquad (j\in\{1,\dots,n\}).
}

\Lem{2}{
Let $f:I\to\R^n$ be an $n$-times continuously differentiable function such that $W_f$ is nonvanishing. Then the components of $f$ form a fundamental system of solutions of the $n$th-order homogeneous linear differential equation 
\Eq{DE+}{
 y^{(n)}
 =\sum_{j=0}^{n-1} \Phi_f^{[j]}y^{(j)}.
}}

\begin{proof}
This equation is equivalent to the following identity
\Eq{*}{
\begin{vmatrix}f^{(n-1)}&\dots&f^{(0)}\end{vmatrix}y^{(n)}
=\sum_{j=0}^{n-1}(-1)^{n-j-1}
\begin{vmatrix}f^{(n)}&\dots &f^{(j+1)}&f^{(j-1)}&\dots&f^{(0)}\end{vmatrix}y^{(j)}.
}
We can now rearrange this equation to obtain
\Eq{*}{
\begin{vmatrix}
y^{(n)}&y^{(n-1)}&\dots&y\\
f_1^{(n)}&f_1^{(n-1)}&\dots&f_1\\
\vdots&\vdots&\ddots&\vdots\\
f_n^{(n)}&f_n^{(n-1)}&\dots&f_n\\
\end{vmatrix}
=0.
}
It is easily seen that if $y\in\{f_1,\dots,f_n\}$, then the determinant vanishes. Therefore, $f_1,\dots,f_n$ are solutions of
\eq{DE+}. Due to the condition that $W_f$ is nonvanishing, the components of $f$ are linearly independent, therefore they 
form a fundamental solution system for \eq{DE+}. 
\end{proof}

\Cor{12}{Let $n,m\in\N$ with $m\geq n$ and let $f:I\to\R^n$ be an $m$-times continuously differentiable function such that $W_f$ is nonvanishing. Define $a:I\to\R^n$ by \eq{aj} and $X_a:I\to\R^{n\times n}$ by \eq{XY}. Then the equality \eq{Wk} holds for $k=(k_1,\dots,k_n)\in\N_0^n$, if $k_i\leq m$ and $\ell_i:=(k_i-n+1)^+$ for $i\in\{1,\dots,n\}$.}

\begin{proof} It follows from the definition of $a$, that it is $(m-n)$-times continuously differentiable. On the other hand, by \lem{2}, we have that $f$ satisfies the $n$-th order homogeneous linear differential equation \eq{DE}. Thus, the statement is a consequence of \thm{1}. 
\end{proof}

We say that two continuous functions $f,g:I\to\R^n$ are \emph{range equivalent}, denoted by $f\sim g$, if there exists a nonsingular $n\times n$-matrix $A$ such that 
\Eq{equ}{
f=Ag.
}

\Thm{3}{Let $f,g:I\to\R^n$ be an $n$-times continuously differentiable functions such that $W_f$ and $W_g$ are nonvanishing. 
Then $f\sim g$ holds if and only if
\Eq{Phi}{
\Phi_f^{[j]}
=\Phi_g^{[j]}
\qquad(j\in\{0,\dots,n-1\}).
}}

\begin{proof}
If $f\sim g$, then there exists a nonsingular $n\times n$-matrix $A$ such that $f=Ag$. The product rule for determinants shows 
that $W_f^k=|A|W_g^k$ for every $k\in\N_0^n$. Using this identity and the definition of $\Phi_f^{[j]}$ and $\Phi_g^{[j]}$, we obtain 
the equalities in \eq{Phi}. 

On the other hand, if the identities \eq{Phi} are valid on $I$, then the $n$th-order homogeneous linear differential equation
\eq{DE+} is equivalent to the following one
\Eq{*}{
 y^{(n)}
=\sum_{j=0}^{n-1} \Phi_g^{[j]}y^{(j)}.
}
Therefore, the ($n$-dimensional) solution spaces of these differential equations are identical, which in view of \lem{2} yields that the components of $f$ are linear combinations of the components of $g$. Thus identity \eq{equ} holds for some nonsingular $n\times n$-matrix $A$.  
\end{proof}


\begin{thebibliography}{1}

\bibitem{Bel27}
E.~T. Bell, \emph{{Partition polynomials}}, Ann. of Math. (2) \textbf{29}
  (1927/28), no.~1-4, 38–46. 

\bibitem{Bel34}
E.~T. Bell, \emph{{Exponential polynomials}}, Ann. of Math. (2) \textbf{35}
  (1934), no.~2, 258–277.

\bibitem{SchRid96}
R.~Schimming and S.~Z. Rida, \emph{{Noncommutative {B}ell polynomials}},
  Internat. J. Algebra Comput. \textbf{6} (1996), no.~5, 635–644.

\bibitem{Tes12}
G.~Teschl, \emph{{Ordinary differential equations and dynamical systems}},
  {Graduate Studies in Mathematics, Vol. 140}, American Mathematical Society,
  2012.
\end{thebibliography}

\end{document}